\documentclass[9pt]{article}
\usepackage{tocloft}

\usepackage{latexsym,amsfonts,euscript,amsthm}

\usepackage{amssymb}
\usepackage{stmaryrd}
\usepackage{mathrsfs,amsmath}
\usepackage{leftidx}
\usepackage{mathabx}

\usepackage[shortlabels]{enumitem}

\usepackage[figuresright]{rotating}
\usepackage[table]{xcolor}
\usepackage{multirow}
\usepackage{booktabs}%tables
\usepackage{tabularx}

\usepackage{tocbibind} %add reference as a section

\usepackage{float}
\usepackage{tikz-cd}
\usetikzlibrary{cd}

\usepackage[pagebackref]{hyperref}

\hypersetup{%
  colorlinks=true,
  linktoc=all,
  citebordercolor=green,
  linkbordercolor=red,
 citecolor=black,
    linkcolor=blue,
  urlcolor=blue
}

\makeatletter
\renewcommand{\@maketitle}{%
  \newpage
  \null
  \vskip 2em%
  \begin{center}%
    \let \footnote \thanks
    {\LARGE \bfseries \@title \par}%  明确加粗
    \vskip 1.5em%
    {\large \textsc{\@author}\par}%
    \vskip 1em%
    {\large \@date \par}%
  \end{center}%
  \par
  \vskip 1.5em
}
\makeatother

\renewcommand{\contentsname}{\small\bfseries Contents}

\newcommand{\centeredtitle}[1]{%
  \begin{center}
    {\small\bfseries #1}%
  \end{center}
  \vspace{1ex}%
}

\makeatletter
\AtBeginDocument{%
  \renewcommand{\tableofcontents}{%
    \centeredtitle{\contentsname}%
    \@starttoc{toc}%
  }%
}
\makeatother

\renewenvironment{abstract}{%
  \centeredtitle{\abstractname}%
  \vspace{-0.6\baselineskip}%  根据实际效果调整负值（如 -1ex, -2ex）
  \quotation
}{%
  \endquotation
}

\newcommand{\ackname}{Acknowledgements}
\makeatletter
\if@titlepage
  \newenvironment{acknowledgement}{%
    \titlepage
    \null\vfil
    \@beginparpenalty\@lowpenalty
    \begin{center}%
      \bfseries \ackname
      \@endparpenalty\@M
    \end{center}}%
  {\par\vfil\null\endtitlepage}
\else
  \newenvironment{acknowledgement}{%
    \if@twocolumn
      \section*{\ackname}%
    \else
      \small
      \begin{center}%
        {\bfseries \ackname\vspace{-0.5em}\vspace{\z@}}%
      \end{center}%
      \quotation
    \fi}
    {\if@twocolumn\else\endquotation\fi}
\fi
\makeatother

\makeatletter % `@' now normal "letter"
\@addtoreset{equation}{section}
\makeatother  % `@' is restored as "non-letter"

\makeatletter
\def\thanks#1{\g@addto@macro\@thanks{\footnotetext{#1}}}
\makeatother

\makeatletter
\newtoks\address@list
\newcommand{\addaddress}[3]{%
  \address@list=\expandafter{\the\address@list
    \textsc{#1, #2.} \par
    \textit{Email address}: \href{mailto:#3}{\textsf{#3}} \par
    \addvspace{\medskipamount}%
  }%
}
\newcommand{\printaddresses}{%
  \AtEndDocument{\bigskip{\footnotesize
    \the\address@list
  }}%
}
\makeatother

\newtheorem{lem}{\bf Lemma}[section]

\newtheorem*{ques}{\bf Question}

\newtheorem*{rmkA}{\bf Remark}
\newtheorem{exa}[lem]{\bf Example}

\newtheorem{mainthm}{Theorem}

\title{An infinite family of counterexamples to a question of Camina
\thanks{
\textbf{2020 MR Subject Classification}\,\, Primary 20D15, 20D60\\
}}

\author{\textsc{Yu Zeng}
}
\addaddress{Yu Zeng}{Department of Mathematics, Suzhou University of Technology, No.99 South Third Ring Road, Changshu, Jiangsu, 215500, China}{yuzeng2004@163.com}
\printaddresses

\date{}

\begin{document}

\maketitle

\begin{abstract}
  A.R. Camina and R.D. Camina posed in \cite{camina06} the following question: 
  Suppose there are two finite groups, one nilpotent and the other non-nilpotent, and the two groups share identical sets of conjugacy class sizes; must the non-nilpotent group possess a non-trivial center?
Recently, W. Zhou \cite{zhou25} gave a negative answer via a subtle and elegant construction of concrete counterexamples.
Nevertheless, his approach relies on the existence of Sophie Germain primes, and thus fails to yield infinitely many counterexamples unconditionally. 
In the present paper, we construct an infinite family of counterexamples to Camina's question.
\end{abstract}

%\tableofcontents

\section{Introduction}

Determining the structural properties of a finite group from its set of conjugacy class sizes
is a classical and interesting problem in group theory.
One origin of this line of study lies in Burnside's $p^\alpha$-lemma,
a fundamental result stating that no finite simple group contains a non-trivial conjugacy class of prime-power size.
Burnside's $p^\alpha$-lemma, together with N. It\^{o}'s celebrated works on finite groups with exactly two conjugacy class sizes, has motivated extensive studies on group characterization via class sizes (for a comprehensive overview, we refer to the survey \cite{camina11}).

For a finite group $G$, let $\mathrm{cs}(G) = \{ |x^G| : x \in G \}$ denote the set of its conjugacy class sizes, where $x^G$ represents the conjugacy class of $x$ in $G$. 
In the literature, instead of the set of conjugacy class sizes, many authors prefer to utilize the notion of the \emph{conjugate type vector} of $G$, which is defined as the tuple of its conjugacy class sizes arranged in ascending order.

Since the set of conjugacy class sizes remains invariant upon taking the direct product of a finite group with any abelian group, class size data alone cannot separate finite simple groups from finite non-simple groups.
However, the well-known Thompson's conjecture on conjugacy class sizes \cite[12.38]{kourovka} addresses
the distinguishability of non-abelian simple groups among centerless finite groups via their class size information.
Another natural question is whether the nilpotency or solvability of a finite group $G$ can be determined solely by $\mathrm{cs}(G)$. 
The answer is negative: A.R. Camina and R.D. Camina \cite{camina06} demonstrated that the set of conjugacy class sizes cannot recognize nilpotency, while G. Navarro \cite{navarro14} established that it also fails to determine solvability.

Much attention has also been paid to conjugacy class sizes of direct products. 
To formulate relevant results clearly,
for two nonempty sets of positive integers $A_1$ and $A_2$, we define their product by
$$ A_1 \cdot A_2 = \{a_1 a_2 : a_1 \in A_1, a_2 \in A_2\}, $$
and write $A_1^m$ for the product of $m$ copies of $A_1$.

It is well known that $\mathrm{cs}(G\times H)=\mathrm{cs}(G)\cdot\mathrm{cs}(H)$ for finite groups $G$ and $H$. 
Since a finite nilpotent group decomposes as the direct product of its Sylow subgroups, its set of class sizes necessarily factors into a product of sets consisting entirely of single-prime powers. 
The converse, however, is false in general.
 A.R. Camina and R.D. Camina \cite{camina06} constructed non-nilpotent groups possessing nilpotent-like class sizes. 
 The smallest such example is a non-nilpotent group $G$ with $\mathrm{cs}(G)=\{1,2,4\}\cdot\{1,5\}$ and $|\mathbf{Z}(G)|=2$.
 In the same paper, they raised the following question.

\begin{ques}[\cite{camina06}, Question 4]\label{q}
  Let $G$ and $H$ be finite groups with $H$ nilpotent. 
  Suppose $\mathrm{cs}(G)=\mathrm{cs}(H)$, but $G$ is not nilpotent. Does $G$ have a non-trivial center?
\end{ques}

Recently, W. Zhou \cite{zhou25} gave a negative answer to this question by constructing a family of centerless groups satisfying $\mathrm{cs}(G) = \{1, p, p^{p-2}\} \cdot \{1, q\}$, where $p$ and $q$ are distinct primes. 
Nevertheless, Zhou's elegant construction relies on a $(p-1)$-dimensional quotient of a permutation module,
imposing the strict constraint $p-1 = 2q$. 
This forces $q$ to be a Sophie Germain prime. 
Since the infinitude of Sophie Germain primes remains unproven, Zhou's method cannot guarantee the existence of infinitely many such counterexamples unconditionally.

Motivated by Zhou's work, we introduce a new construction based on truncated polynomial rings, which successfully removes the restriction on $p$ and $q$. 
Combined with Dirichlet's theorem on arithmetic progressions, our approach yields an unconditional construction of infinitely many counterexamples to Camina's problem. 
Our main result is stated as follows.

\begin{mainthm}\label{thmA}
Let $q$ be an arbitrary prime, and let $p$ be any prime such that $p \equiv 1 \pmod q$ and $p \ge q+2$. There exists a finite $\{ p,q \}$-group $G$ with trivial center, such that $\mathrm{cs}(G) = \{1, p, p^{q+1}\} \cdot  \{1, q\}$. 
In particular, $G$ is centerless, but it shares its set of conjugacy class sizes with a finite nilpotent group.  
\end{mainthm}

The above result leads to several further observations and extensions.

\begin{rmkA}
  {\rm 
  By Dirichlet's theorem on arithmetic progressions \cite[II, Theorem 15]{hardy08}, for every prime $q$, there exist infinitely many primes $p \equiv 1 \pmod q$, which guarantees that Theorem \ref{thmA} yields an infinite family of counterexamples. 
  Moreover, the flexibility of our construction allows for several extended structural properties. 
  It enables the creation of centerless groups $G$ with an arbitrarily large number of distinct prime divisors $|\pi(G)|$, as well as groups satisfying $\mathrm{cs}(G)=\mathrm{cs}(P)\cdot \mathrm{cs}(Q)$ for a finite $p$-group $P$ and finite $q$-group $Q$ (where $p$ and $q$ are distinct primes) with arbitrarily large 
  $|\mathrm{cs}(P)|$ and $|\mathrm{cs}(Q)|$.
  Indeed, it is possible to construct groups satisfying both conditions simultaneously. 
  We explicitly demonstrate these properties through the following constructions.
 \begin{description}
  \item[(1)] For every positive integer $n$, we can choose $2n$ distinct primes $p_1, q_1, \dots, p_n, q_n$ such that $p_i \equiv 1 \pmod{q_i}$ and $p_i \geq q_i+2$ for all $1 \leq i \leq n$. 
  Let $G_i$ be the finite $\{ p_i, q_i \}$-group in Theorem \ref{thmA}, and define $G=G_1\times \cdots \times G_n$. 
  Then $G$ is a finite group with a trivial center such that $\mathrm{cs}(G)=\mathrm{cs}(H)$ for some finite nilpotent group $H$ by Lemma \ref{lem:cossey}. 
   Since the chosen primes are mutually distinct, we have $|\pi(G)|=2n$, demonstrating that the number of prime divisors of $|G|$ can be arbitrarily large.

  \item[(2)] Let $G_1$ be the finite $\{ p,q \}$-group in Theorem \ref{thmA}. For every positive integer $m$, let $G_m$ be the direct product of $m$ copies of $G_1$. Then $G_m$ has a trivial center, and its set of conjugacy class sizes is given by
  \[ \mathrm{cs}(G_m) = \mathrm{cs}(G_1)^m = \{1, p, p^{q+1}\}^m \cdot \{1, q\}^m. \]
  Evaluating these sets, we observe that $\{1, q\}^m = \{1, q, q^2, \dots, q^m\}$ contains exactly $m+1$ distinct sizes. Similarly, the set $\{1, p, p^{q+1}\}^m$ contains at least $m+1$ distinct powers of $p$ (for instance, the subset $\{1, p, p^2, \dots, p^m\}$). Therefore, by Lemma \ref{lem:cossey}, we can write $\mathrm{cs}(G_m) = \mathrm{cs}(P) \cdot \mathrm{cs}(Q)$, where $P$ is a finite $p$-group and $Q$ is a finite $q$-group, such that both $|\mathrm{cs}(P)|$ and $|\mathrm{cs}(Q)|$ grow unboundedly as $m \to \infty$. 
 \end{description}
  }
\end{rmkA}

The organization of the paper is as follows. 
Section 2 collects necessary preliminary results. 
In Section 3, we construct examples and verify that they satisfy Theorem \ref{thmA}.

Throughout the rest of this paper, we only consider finite groups.

\section{Preliminaries}

We rely on \cite{huppert67} as our primary reference for standard results in finite group theory and representation theory. 
Throughout the paper, $p$ and $q$ always denote distinct primes, and $\mathbb{F}_p$ denotes the finite field of $p$ elements. 
Any remaining notation will be introduced as needed.

A key tool for matching our constructed set $\mathrm{cs}(G)$ with that of a finite nilpotent group is the following classical theorem due to J. Cossey and T. Hawkes.

\begin{lem}\label{lem:cossey}
Every finite set of $p$-powers containing $1$ is the set of conjugacy class sizes of some finite $p$-group.
\end{lem}
\begin{proof}
  This is \cite[Theorem]{cossey00}.
\end{proof}

Let $G$ be a finite group, and let $V$ be a $G$-module which we always write additively in the following. 
For every $g\in G$ and $v\in V$, the image of $v$ under the action of $g$ is denoted by $v.g$.

\begin{lem}\label{lem: frob}
  Let $G$ be a Frobenius group with kernel $A$ and complement $B$, and let $V$ be a faithful $G$-module over the prime field $\mathbb{F}_p$ with no non-trivial $G$-fixed points.
Suppose $A=\langle a\rangle\cong \mathsf{C}_p$ and $B=\langle b\rangle\cong \mathsf{C}_q$ for some prime $q$. Then
$$\left|\bigcup_{i=0}^{p-1} \mathbf{C}_{V}(b^{a^i}) \cup \mathbf{C}_{V}(a) \right|=|\mathbf{C}_{V}(a)|+p |\mathbf{C}_{V}(b)|-p.$$
Moreover, either the set $\bigcup_{i=0}^{p-1}\mathbf{C}_V(b^{a^i})\cup \mathbf{C}_V(a)$ is a proper subset of $V$, or $|V|=p^{2}$.
In the latter case, by suitably identifying $V$ with $\mathbb{F}_p\oplus\mathbb{F}_p$ (the space of column vectors),
 the group $G$ embeds into $\mathrm{GL}_2(p)$ (acting via left multiplication on $\mathbb{F}_p\oplus\mathbb{F}_p$) as
\[
G=\left\{
\begin{pmatrix}
e & d\\
0 & 1
\end{pmatrix}\;:\; e\in E,\ d\in\mathbb{F}_p
\right\},
\]
where $E\leq \mathbb{F}_p^\times$ is of order $q$.
\end{lem}
\begin{proof}
 Let $X=\{ a, b^{a^i}: 0\leq i\leq p-1 \}$. For each pair of elements $x,y\in X$, we have $G=\langle x,y\rangle$.
Since $|\mathbf{C}_{V}(G)|=1$, it follows that $|\mathbf{C}_{V}(x)\cap \mathbf{C}_{V}(y)|=1$.
Because these $p+1$ subgroups ($\mathbf{C}_{V}(a)$ and the $p$ conjugates of $\mathbf{C}_{V}(b)$) intersect trivially pairwise, the size of their union is simply the sum of their non-trivial elements plus one:
\[
\left|\bigcup_{i=0}^{p-1} \mathbf{C}_{V}(b^{a^i}) \cup \mathbf{C}_{V}(a)\right| = 1 + (|\mathbf{C}_{V}(a)| - 1) + p(|\mathbf{C}_{V}(b)| - 1) = |\mathbf{C}_{V}(a)| + p|\mathbf{C}_{V}(b)| - p.
\]

Assume this union covers $V$. Then we have:
\[
|V| - |\mathbf{C}_{V}(a)| = p(|\mathbf{C}_{V}(b)| - 1).
\]
Since $G$ acts faithfully on $V$, $1<|\mathbf{C}_{V}(b)| < |V|$. 
 The equation above then forces $|\mathbf{C}_{V}(a)|=p$ and $|V| = p|\mathbf{C}_{V}(b)|$.
Because the $p$ subgroups $\mathbf{C}_{V}(b^{a^i})$ intersect trivially pairwise and each has dimension $\dim_{\mathbb{F}_p}(V) - 1$, we must have $|\mathbf{C}_{V}(b)|=p$.
Since $|\mathbf{C}_{V}(a)\cap \mathbf{C}_{V}(b)|=1$, we deduce that $V=\mathbf{C}_{V}(a)\oplus \mathbf{C}_{V}(b)$ with $|V|=p^{2}$.

Let $v \in \mathbf{C}_{V}(a)$ and $w \in \mathbf{C}_{V}(b)$ be non-trivial elements.
Since $\mathbf{C}_{V}(a)$ is $G$-invariant and intersects $\mathbf{C}_{V}(b)$ trivially, $\langle b \rangle\cong \mathsf{C}_{q}$ acts fixed-point-freely on $\mathbf{C}_{V}(a)$, so $v.b=ev$ for some $e\in \mathbb{F}_p^\times$ of order $q$.
Furthermore, since $a$ acts trivially on both $\mathbf{C}_{V}(a)$ and $V/\mathbf{C}_{V}(a)$, we have $v.a=v$ and $w.a=dv+w$ for some $d\in \mathbb{F}_p$.
Identifying $v$ and $w$ with the (basis) column vectors $(1,0)^T$ and $(0,1)^T$ in $\mathbb{F}_p\oplus \mathbb{F}_p$, the final statement easily follows.
\end{proof}

\section{Construction and proofs}

In this section,
we establish Theorem \ref{thmA} by constructing $G = V \rtimes (A \rtimes B)$ via a Frobenius group $A \rtimes B$ acting on a quotient ring $V$ of the one-variable polynomial ring $\mathbb{F}_p[X]$ over a prime field $\mathbb{F}_p$.

\begin{exa}\label{exa}
For every prime $q$, choose a prime $p \equiv 1 \pmod q$ with $p \geq q+2$, and let 
$r<p$ be a positive integer of multiplicative order $q$ modulo $p$.
We construct the group $G = V \rtimes (A \rtimes B)$ as follows.
\begin{description}
  \item[(1)] First, let $A = \langle a \rangle \cong \mathsf{C}_p$ and $B = \langle b \rangle \cong \mathsf{C}_q$, where $B$ acts on $A$ via $b^{-1} a b = a^r$. Thus, $A \rtimes B$ is a Frobenius group with kernel $A$.
Next, let $V=\mathbb{F}_p[x] =\mathbb{F}_p[X] / (X^{q+2})$, where $(X^{q+2})$ denotes the ideal generated by $X^{q+2}$ and $x:=X+(X^{q+2})$. Viewed as an additive group, $V$ is an elementary abelian $p$-group of order $p^{q+2}$.
  \item[(2)] The linear action of $A \rtimes B$ on $V$ is defined for every polynomial $f(x) \in V$ by
\begin{center}
  $f(x).a =  f(x)(1+x)$,~~~~~~~~~~$f(x). b=f((1+x)^r - 1)$.
\end{center}
\end{description}
\end{exa}

We first show that the group defined in Example \ref{exa} is well-defined.

\begin{lem}\label{lem: well define}
The group $G$ defined in Example {\rm\ref{exa}} is well-defined. Moreover, $|G|=p^{q+3}q$ and $V=\mathbf{C}_{G}(V)$.
\end{lem}
\begin{proof}
We only need to verify that the defined actions preserve all defining relations of the group $AB$. 
Since $p \ge q+2$, $x^p = x^{q+2}=0 $. Thus, 
\[
f(x). a^p = f(x)(1+x)^p = f(x)(1+x^p) =f(x),
\]
confirming that $a$ is a linear transformation of order $p$ on $V$.
Note that, for every positive integer $k$,
\[
f(x).b^k =  f((1+x)^{r^k} - 1),
\]
and so $b^{q}$ acts on $f(x)$ via substituting $x$ with $(1+x)^{r^q} - 1$.
Since $r^q \equiv 1 \pmod p$, we can write $r^q = 1 + cp$ for some integer $c$. 
Thus $(1+x)^{1+cp}-1 =  (1+x)(1+x^p)^c-1 =  x$, which shows that $b$ is a linear transformation of order $q$ on $V$.
Finally, we verify the semidirect product relation.
Note that
\[ f(x).b^{-1}ab = f((1+x)^{r^{q-1}}-1).ab = f((1+x)^{r^{q-1}}-1)(1+x).b =   f(x) ((1+x)^r - 1 + 1) = f(x)(1+x)^r = f(x). a^{r}.\]
Hence $b^{-1}ab = a^r$ on $V$, and $G$ is well-defined. 
Moreover, $|G|=|V||A||B| = p^{q+3}q$ and $V=\mathbf{C}_{G}(V)$.
\end{proof}

We are now ready to prove Theorem \ref{thmA}.

\begin{proof}[Proof of Theorem \ref{thmA}]
 Let $G$ be the finite group defined in Example \ref{exa}.
  We proceed with the proof by verifying the claims in steps.

\textbf{Step 1.} $|\mathbf{C}_{V}(a)|=p$, and $\mathbf{Z}(G)=\mathbf{C}_{V}(a)\cap \mathbf{C}_{V}(b)$ has order 1.

Let $f(x)\in \mathbf{C}_{V}(a)$. Since $f(x)(1+x)= f(x) $, we must have $xf(x)= 0$.
Thus, $\mathbf{C}_{V}(a)=\{ c x^{q+1} : c\in \mathbb{F}_p \}$.
In particular, $|\mathbf{C}_{V}(a)|=p$.
Since $V=\mathbf{C}_{G}(V)$ by Lemma \ref{lem: well define}, 
 we have $\mathbf{Z}(G)\leq \mathbf{C}_{G}(V)=V$.
Moreover,
it is routine to verify that $\mathbf{Z}(G)=\mathbf{C}_{V}(a)\cap \mathbf{C}_{V}(b)$.
Observe the action of $b$ on the basis element $x^{q+1}$ of $\mathbf{C}_{V}(a)$:
  \[
  x^{q+1}.b= ((1+x)^{r}-1)^{q+1}=(r x+\cdots )^{q+1}=r^{q+1}x^{q+1}.
  \]
Since $r$ has multiplicative order $q\geq 2$ modulo $p$, it follows that $r^{q+1} \equiv r \not\equiv 1 \pmod{p}$.
Therefore, $b$ acts fixed-point-freely on $\mathbf{C}_{V}(a)$, implying $|\mathbf{Z}(G)| = |\mathbf{C}_{V}(a)\cap \mathbf{C}_{V}(b)|=1$.

\textbf{Step 2.} $\{ |v^{G}| : v\in V \}=\{ 1, p, q, pq\}$.

  If $v$ is the trivial element of $V$, then $|v^{G}|=1$.
  Now let $v$ be a non-trivial element of $V$.
Since $V$ is abelian, $V \leq \mathbf{C}_{G}(v)$, so the centralizer splits as $\mathbf{C}_{G}(v) = V \rtimes \mathbf{C}_{AB}(v)$.
Because $|AB|=pq$ with $p$ and $q$ distinct primes, the index $[AB : \mathbf{C}_{AB}(v)]$ must be $1, p, q,$ or $pq$.
Consequently, the size $|v^G| = [AB : \mathbf{C}_{AB}(v)]$ belongs to $\{p, q, pq\}$.
For non-trivial elements $u\in \mathbf{C}_{V}(a)$ and $w\in \mathbf{C}_{V}(b)$, the fact that $\mathbf{Z}(G)=1$ ensures neither is fully centralized by $AB$. Thus, $|u^{G}|=q$ and $|w^{G}|=p$.
Finally, since $|V|=p^{q+2}>p^{2}$, Lemma \ref{lem: frob} guarantees the existence of an element $x\in V- (\bigcup_{i=0}^{p-1}\mathbf{C}_{V}(b^{a^{i}})\cup \mathbf{C}_{V}(a))$.
For this $x$, the centralizer in $AB$ is trivial, so $\mathbf{C}_{G}(x)=V$, which concludes that $|x^{G}|=pq$.

\textbf{Step 3.} If $v\in V$ and $i\not\equiv 0 \pmod{p}$, then $|(va^{-i})^{G}|=p^{q+1}q$.    

An element $w a^{-j} b^k \in G$ (where $w \in V$) centralizes $v a^{-i}$ if and only if $k\equiv 0 \pmod{q}$ (since $B$ acts fixed-point-freely on $A$) and the vector equation $w.(a^i-1) = v.(a^j-1)$ holds in $V$. 

%Without loss of generality, we may assume that $0<i<p$.
Since $i \not\equiv 0 \pmod p$, the linear transformation $a^i-1$ corresponds to multiplication by $(1+x)^i-1 = x(i + \cdots)$. In the quotient ring $V$, this multiplication maps $V$ precisely onto the ideal $(x)$ generated by $x$, which means the kernel of the linear transformation $a^{i}-1$ is 1-dimensional. Thus, $|\ker(a^i-1)| = p$.
For every choice of $j$, the linear transformation $a^j-1$ corresponds to multiplication by $(1+x)^j-1$, which also maps $V$ into the ideal $(x)$. Consequently, the vector $v.(a^j-1)$ always lies in the image of $a^i-1$. 
Viewing $w.(a^i-1) = v.(a^j-1)$ as an equation for unknown $w$ in the space $V$,
we deduce that the equation has exactly $|\ker(a^i-1)| = p$ solutions.
% This guarantees that the equation $w.(a^i-1) = v.(a^j-1)$ is consistently solvable for $w$??, yielding exactly $|\ker(a^i-1)| = p$ solutions. 
Since there are exactly $p$ choices for $j$, we have $|\mathbf{C}_G(v a^{-i})| = p^2$.
As a consequence, $|(v a^{-i})^{G}|=[G:\mathbf{C}_{G}(va^{-i})] = p^{q+3}q / p^2 = p^{q+1}q$.

\textbf{Step 4.} If $v\in V$ and $j\not\equiv 0 \pmod{q}$, then $|(va^{i}b^{-j})^{G}|=p^{q+1}$.    

Since $A B$ is a Frobenius group with complement $B=\langle b^{j}\rangle$, by \cite[V, Satz 8.9]{huppert67}
the commutator map $z \mapsto [z, b^{j}]$ is a bijection on $A$.
In particular, there exists some $y\in A$ such that $[y,b^{j}]=a^{i}$,
meaning $a^i b^{-j}=(b^{-j})^{y}$.
%any element $a^i b^{-j}$ is $A$-conjugate to $b^{-j}$.
%for some integer $k\not\equiv 0 \pmod{q}$. 
Consequently, $v a^i b^{-j}$ is $G$-conjugate to $w b^{-j}$ for some $w\in V$.
%for some $w \in V$ and some integer $k\not\equiv 0 \pmod{q}$.
 Thus, it suffices to compute the size of $(v b^{-j})^G$.
An element $w a^k b^{-n} \in G$ (where $w \in V$) centralizes $v b^{-j}$ if and only if $k\equiv 0\pmod{p}$ (since $B$ acts fixed-point-freely on $A$) and the vector equation $w.(b^j-1) = v.(b^n-1)$ holds in $V$. 

To solve this equation, consider the action of $b$ on the basis $\{ x^m : 0 \leq m \leq q+1 \}$ of $V$. We have
\[
x^m. b = ((1+x)^r-1)^m = r^m x^m + \cdots 
\]
Since the order of $b$ is coprime to $p$,
it follows that $b$ acts as a diagonalizable linear transformation on $V$ with eigenvalues $1, r, r^2, \dots, r^{q+1}$ (viewing as elements in $\mathbb{F}_p$). Since $r$ has multiplicative order $q$ modulo $p$, the eigenvalue $1$ occurs exactly when $m=0$ or $m=q$. Thus, the 1-eigenspace $\mathbf{C}_V(b)$ of $b$ has dimension 2, giving $|\mathbf{C}_V(b)| = p^2$. 
Again, as $\langle b\rangle$ acts coprimely on $V$, we decompose $V$ into the direct sum $V = \mathbf{C}_V(b) \oplus [V, b]$.

The linear transformations $b^j-1$ and $b^n-1$ annihilate $\mathbf{C}_V(b)$, which ensures that both the vectors $w.(b^{j}-1)$ and $v.(b^n-1)$ lie in $[V, b]$. 
Furthermore, since $j \not\equiv 0 \pmod q$, the linear transformation $b^j$ has no eigenvalue equal to 1 on $[V, b]$, making $b^j-1$ invertible on this subspace. Thus, for each of the $q$ choices of $n$, the equation $w.(b^j-1) = v.(b^n-1)$ is solvable and yields exactly $|\ker(b^j-1)| = |\mathbf{C}_V(b)| = p^2$ solutions for $w$.
This gives $|\mathbf{C}_G(v b^{-j})| = p^2 q$, and therefore $|(v b^{-j})^G| = [G : \mathbf{C}_G(v b^{-j})] = p^{q+3}q / p^2q = p^{q+1}$.

\textbf{Step 5. Final Conclusion.} 

Combining the above steps, we deduce that $G$ is centerless, and its set of conjugacy class sizes is given by 
\[ \mathrm{cs}(G) = \{1, p, q, pq, p^{q+1}, p^{q+1}q\} = \{1, p, p^{q+1}\} \cdot \{1, q\}. \]
Thus, by Lemma \ref{lem:cossey}, there exists a finite nilpotent group $H$ such that $\mathrm{cs}(G) = \mathrm{cs}(H)$.
\end{proof}

\begin{acknowledgement}
	%The author would like to thank Professor Guohua Qian for drawing my attention to this problem.
   The author is grateful to the referee for her/his valuable comments.
\end{acknowledgement}

\end{document}